# On universal and epi-universal locally nilpotent groups [*]

R. Göbel, S. Shelah and S.L. Wallutis [**]

**Abstract.** In this paper we mainly consider the class $\mathcal{LN}$ of all locally nilpotent groups. Using similar arguments as in [GrS] we first show that there is no universal group in $\mathcal{LN}_\lambda$ if $\lambda$ is a cardinal such that $\lambda = \lambda^{\aleph_0}$; here we call a group $G$ *universal* (in $\mathcal{LN}_\lambda$) if any group $H \in \mathcal{LN}_\lambda$ can be embedded into $G$ where $\mathcal{LN}_\lambda$ denotes the class of all locally nilpotent groups of cardinality at most $\lambda$.

However, our main interest is the construction of torsion-free epi-universal groups in $\mathcal{LN}_\lambda$, where $G \in \mathcal{LN}_\lambda$ is said to be *epi-universal* if any group $H \in \mathcal{LN}_\lambda$ is an epimorphic image of $G$. Thus we give an affirmative answer to a question by Plotkin (see [KN, 3.47]). To prove the torsion-freeness of the constructed locally nilpotent group we adjust the well-known commutator collecting process due to P. Hall to our situation.

Finally, we briefly discuss how to use the same methods as for the class $\mathcal{LN}$ for other canonical classes of groups to construct epi-universal objects.

## Introduction

Many categories $\mathcal{K}$ in algebra or geometry share the property that, given a cardinal $\lambda$, one can find a *universal* object $G \in \mathcal{K}$, that is, any object in $\mathcal{K}$ of cardinality $\leq \lambda$ embeds into $G$.

In this paper we shall restrict our attention to classes of groups. If the class of groups is sufficiently homogeneous and if there are many embeddings, then universal objects do exist. For example, in the case of abelian groups the universal objects are obviously certain divisible groups. However, if we restrict the number of embeddings, replacing monomorphisms by pure injections, then it can be shown that universal

[*]GbSh 742 in Shelah's list of publications
[**]supported by Deutsche Forschungsgemeinschaft

objects often do not appear; see the recent works by Kojman, Shelah [KS] and Shelah [S].

On the other hand, if the class of groups is less homogeneous, then it seems hopeless to search for universal objects. A good example is the class of locally finite groups, as shown by Grossberg, Shelah [GrS]. Modifying the proof in [GrS] we note that the class of locally nilpotent groups has no universal objects (see Theorem 1.3).

It is thus natural to consider the dual notion. This is to say that an object $G$ of a given class $\mathcal{K}$ of groups with $|G| = \lambda$ is *epi-universal* in the subclass $\mathcal{K}_\lambda$ of $\mathcal{K}$ consisting of all those elements of cardinality at most $\lambda$, if any group in $\mathcal{K}_\lambda$ is an epimorphic image of $G$. Clearly, all group varieties $\mathcal{V}$ have epi-universal objects, namely their $\mathcal{V}$-free groups. However, the class $\mathcal{LN}$ of all locally nilpotent groups does not form a variety. Recall, that a group $G$ is said to be *locally nilpotent*, if any finitely generated subgroup $H \subseteq G$ is nilpotent, i.e. the lower central series $H^0 = H \supseteq H^1 = [H, H] \supseteq \ldots \supseteq H^{k+1} = [H^k, H] \ldots$ is finite.

It is the aim of this paper to construct particular "pseudo-free" locally nilpotent groups, which can be shown to be torsion-free by modification of Hall's commutator collecting process. Moreover, any locally nilpotent group is the epimorphic image of a suitable pseudo-free locally nilpotent group, which answers a problem by Plotkin (see [KN, Problem 3.47]).

For certain cardinals $\lambda$, namely cardinals satisfying $\lambda = \lambda^{<\lambda}$, we show that there is a pseudo-free locally nilpotent group of cardinality $\lambda$ which is also epi-universal in the class $\mathcal{LN}_\lambda$ of all locally nilpotent groups of cardinality $\leq \lambda$. Note that, assuming GCH, all regular cardinals satisfy the above condition.

Finally, an analysis of the proofs shows that the existence of epi-universal objects follows similarly for the class $\mathcal{LV}$ consisting of all $\mathcal{V}$-locally groups for any union $\mathcal{V}$ of ascending chains $\{\mathcal{V}_n \mid n < \omega\}$ of group varieties $\mathcal{V}_n$, where we call a group $G$ $\mathcal{V}$-*locally* if all its finitely generated subgroups are elements of $\mathcal{V}$ (cf. [R]). Clearly, $\mathcal{LN}$ is such a class taking $\mathcal{V}_n$ to be the variety of all nilpotent groups of class $\leq n$. Another example are the locally solvable groups with $\mathcal{V}_n$ being all solvable groups of length at most $n$. Also, the class of all torsion groups can obviously be described in this way and thus the existence of epi-universal torsion groups follows. Note, however, that there are no universal torsion groups for $\lambda = \lambda^{\aleph_0}$ (see Proposition 1.4).

We use standard notations from group theory with the exception that our maps act on the right.

## §1 Non-existence of universal locally nilpotent groups

In this section we show that, for certain cardinals $\lambda$, there is no universal member in the class $\mathcal{LN}_\lambda$ of all locally nilpotent groups of cardinality at most $\lambda$.

The existence of universal objects in any class is strongly related to the amalgamation property (see [BS, Ch. 11, §3]). Thus, in order to show that there is no universal element, we prove the failure of the amalgamation property by giving an explicit example. Recall, a class $\mathcal{K}$ satisfies the amalgamation property if, for any $A, B_0, B_1 \in \mathcal{K}$ with embeddings $f_0 \colon A \hookrightarrow B_0$, $f_1 \colon A \hookrightarrow B_1$, there are $G \in \mathcal{K}$ and embeddings $h_0 \colon B_0 \hookrightarrow G$, $h_1 \colon B_1 \hookrightarrow G$ such that $f_0 h_0 = f_1 h_1$.

We begin with constructing suitable groups.

**Construction 1.1** Let $A$ be the elementary abelian 2-group generated by the independent set $\{x_n \mid n < \omega\}$. We define mappings $g_0, g_1$ on $A$ by

$$x_n g_0 = \begin{cases} x_n x_{n+1} & ; \text{ for } n \text{ even} \\ x_n & ; \text{ for } n \text{ odd} \end{cases}$$

and

$$x_n g_1 = \begin{cases} x_n & ; \text{ for } n \text{ even} \\ x_n x_{n+1} & ; \text{ for } n \text{ odd}. \end{cases}$$

Clearly, the mappings extend to automorphisms $g_0, g_1$ of $A$ with $g_0^2 = g_1^2 = 1$.

We now define subgroups $B_0$, $B_1$ of the holomorph of $A$ as follows (see [R, p.37] for the definition of holomorphs):

Let $A^*$ denote the right regular permutations of $A$ and put $B_0 = \langle A^*, g_0 \rangle$ and $B_1 = \langle A^*, g_1 \rangle$. Obviously, the following relations are satisfied:

$$(x_n^*)^{g_0} = \begin{cases} x_n^* \, x_{n+1}^* & ; \text{ for } n \text{ even} \\ x_n^* & ; \text{ for } n \text{ odd} \end{cases} \quad (1)$$

and

$$(x_n^*)^{g_1} = \begin{cases} x_n^* & ; \text{ for } n \text{ even} \\ x_n^* \, x_{n+1}^* & ; \text{ for } n \text{ odd} \end{cases} \quad (2)$$

where $(x_n^*)^{g_i} = g_i x_n^* g_i$ denotes the conjugate.

Note, $g_0$ and $g_1$ are now inner automorphisms of $B_0$ and $B_1$, respectively. For simplification we identify $A^*$ with $A$ and write $B_0 = \langle A, y_0 \rangle$, $B_1 = \langle A, y_1 \rangle$ where $y_0, y_1$ are elements of order 2 satisfying the corresponding relations (1) and (2). □

We now show that these groups are, indeed, suitable for proving the failure of the amalgamation property.

**Proposition 1.2** *The groups $B_0, B_1$ as defined in Construction 1.1 are locally finite 2-groups; in particular, they are locally nilpotent. Moreover, the maps $\mathrm{id}_A \colon A \hookrightarrow B_0$, $\mathrm{id}_A \colon A \hookrightarrow B_1$ do not satisfy the amalgamation property in $\mathcal{LN}$.*

**Proof.** First we show that $B_0$ is a locally finite 2-group. Let $U$ be any finitely generated subgroup of $B_0$. Then we can choose an odd integer $k$ such that $U \subseteq \langle y_0, x_0, x_1, \ldots, x_k \rangle =: V$. Let $v \in V$. Then either $v \in A_k = \langle x_n \mid n \leq k \rangle$ or $v = y_0 w$ for some $w \in A_k$ where the latter follows from $x_n^{y_0} \in A_k$ for all $n \leq k$. Clearly, the order of $v$ is $2^l$ for some $l$ and thus $V$ is a finite 2-group. Therefore $B_0$, and similarly $B_1$, is a locally finite 2-group.

It remains to show that for any locally nilpotent group $G$ there are no embeddings $h_0 \colon B_0 \hookrightarrow G$ and $h_1 \colon B_1 \hookrightarrow G$ such that $h_0 \restriction A = h_1 \restriction A$. Suppose, for contradiction, that there are such embeddings and let $H$ be the subgroup of $G$ generated by the finite set $\{x_0 h_0 = x_0 h_1, \, z_0, \, z_1\}$ where $z_0 = y_0 h_0$, $z_1 = y_1 h_1$. Now,

$$x_0 h_0 \in H = H^0,$$
$$x_1 h_0 = ([x_0, y_0])h_0 = [x_0 h_0, z_0] \in H^1 = [H, H],$$
$$x_2 h_0 = x_2 h_1 = [x_1 h_1, z_1] \in H^2 = [H^1, H] \ldots$$

and so on, i.e.
$$x_n h_0 \in H^n \text{ for all } n < \omega.$$

Therefore $H^n \neq 1$ since $x_n h_0 \neq 1$ by the injectivity of $h_0$ and so $H$ is not nilpotent contradicting the local nilpotency of $G$. $\square$

We are now ready to prove the main result of this section:

**Theorem 1.3** *Let $\lambda$ be a cardinal with $\lambda = \lambda^{\aleph_0}$.*

*Then there is no universal group in $\mathcal{LN}_\lambda$.*

**Proof.** Let $G$ be an arbitrary locally nilpotent group of cardinality $\lambda$. We shall construct a group $H \in \mathcal{LN}_\lambda$ which cannot be embedded into $G$.

Let $T$ be the tree given by $T = {}^{\omega >}\lambda = \{\tau \colon n \to \lambda \mid n < \omega\}$. Moreover, let $G_0$ be the elementary abelian 2-group defined by

$$G_0 := \langle x_\tau \mid \tau \in T \rangle \text{ with } x_\tau^2 = 1 \ (\tau \in T).$$

For any $v \colon \omega \to \lambda$ we put $K_v := \langle x_\tau \mid \tau \in v \rangle \subseteq G_0$ where we identify $v$ with $\{v \upharpoonright n \mid n < \omega\} \subseteq T$. Note, the group $K_v$ is isomorphic to the group $A$ as in the Construction 1.1. We first show:

$$\begin{array}{l} \text{There exists a family } \{h_v \colon K_v \hookrightarrow G \mid v \in {}^\omega\lambda\} \\ \text{of embeddings such that, for any embedding} \\ h \colon G_0 \hookrightarrow G, \text{ there is } v \in {}^\omega\lambda \text{ with } h_v \subseteq h. \end{array} \quad (+)$$

To prove (+) we inductively define monomorphisms $h_\tau \colon \langle x_{\tau \upharpoonright n} \mid n < \operatorname{dom} \tau \rangle \hookrightarrow G$ for any $\tau \in T$. For $\operatorname{dom} \tau = 0$ we put $h_\tau = 0$. Suppose $h_\tau$ is already defined for some $\tau \in T$ with $\operatorname{dom} \tau = k$. We enumerate the set of all embeddings

$$\varphi \colon \langle x_{\tau \upharpoonright n} \mid n \leq k \rangle \hookrightarrow G \text{ with } h_\tau \subseteq \varphi$$

by
$$\{\varphi_\alpha \mid \alpha < \lambda\}.$$

For any $\alpha < \lambda$ we define
$$h_{\tau^\frown \langle \alpha \rangle} \colon \langle x_{\tau \upharpoonright n} \mid n \leq k \rangle \hookrightarrow G$$

by
$$h_{\tau\hat{}\langle\alpha\rangle} \restriction \langle x_{\tau \restriction n} \mid n < k \rangle = h_\tau$$
and
$$x_\tau h_{\tau\hat{}\langle\alpha\rangle} = x_\tau \varphi_\alpha \quad (\tau = \tau \restriction k)$$

where $\tau\hat{}\langle\alpha\rangle \in T$ denotes the mapping $\tau' \colon k+1 \to \lambda$ with $\tau \subseteq \tau'$ and $k\tau' = \alpha$. Since $h_\tau \subseteq h_\sigma$ for $\tau \subseteq \sigma$ (in $T$) we may now define $h_v = \bigcup_{n<\omega} h_{v\restriction n}$ for all $v \in {}^\omega\lambda$. Next we show that the $h_v$'s really satisfy (+). To do so let $h \colon G_0 \hookrightarrow G$ be any embedding. We obtain $v = \{v_n = v \restriction n \mid n < \omega\}$ as follows: For $n = 0$ we have no choice; $v_0 = \emptyset$ (empty mapping). Suppose we have found $v_k$ for some $k < \omega$ with $h_{v_k} \subseteq h$. By the definition of the $h_\tau$'s there exists an ordinal $\alpha < \lambda$ such that $\varphi_\alpha = h \restriction \langle x_{v \restriction n} \mid n \leq k\rangle = h_{(v_k)\hat{}\langle\alpha\rangle}$. We put $v_{k+1} = v_k\hat{}\langle\alpha\rangle$. Obviously, we thus have $h_v = \bigcup_{n<\omega} h_{v_n} \subseteq h$ as required, i.e. (+) is proved.

We now want to apply Proposition 1.2 to the countable subgroups $K_v$ ($v \in {}^\omega\lambda$) of $G_0$. For any $v \in {}^\omega\lambda$ and $\varepsilon = 0, 1$ we define elements $y_\varepsilon^v$ as in Construction 1.1 by
$$(y_\varepsilon^v)^2 = 1$$
and by the conjugates
$$x_{v \restriction k}^{y_\varepsilon^v} = y_\varepsilon^v \, x_{v \restriction k} \, y_\varepsilon^v = x_{v \restriction k} \, x_{v \restriction (k+1)}$$
for $k$ even, $\varepsilon = 0$ and for $k$ odd, $\varepsilon = 1$; otherwise we put
$$x_\tau^{y_\varepsilon^v} = x_\tau.$$

By Proposition 1.2 there are $\varepsilon_v \in \{0,1\}$ ($v \in {}^\omega\lambda$) such that $h_v$ cannot be extended to an embedding from $\langle K_v, y_{\varepsilon_v}^v\rangle$ into $G$. For simplicity we write $y_v = y_{\varepsilon_v}^v$.

Now let $H$ be defined by
$$H := \langle G_0, \, y_v \mid v \in {}^\omega\lambda\rangle.$$

Clearly, $|H| = \lambda$ since $|{}^\omega\lambda| = \lambda^{\aleph_0} = \lambda = |T|$.

Next we show that $H$ is a locally finite 2-group and hence $H \in \mathcal{LN}_\lambda$.

A finitely generated subgroup $U$ of $H$ looks like
$$U = \langle x_\tau \, (\tau \in E), \, y_v \, (v \in F)\rangle$$
for some finite sets $E \subseteq T, F \subseteq {}^\omega\lambda$.

Let $k < \omega$ be such that $v \restriction k \neq w \restriction k$ for all $v \neq w \in F$. W.l.o.g. we may assume that $v \restriction n \in E$ for all $n \leq k$ and all $v \in F$. Moreover, we may assume that, for any $v \in F$, the element $\tau$ of maximal domain ($\geq k$) in $E \cap v$ satisfies $x_\tau^{y_v} = x_\tau$ and $\tau \restriction n \in E$ for all $n \leq \operatorname{dom} \tau$. Therefore we have that $x_\sigma^{y_v} \in U_E = \langle x_\tau \mid \tau \in E \rangle$ for any $\sigma \in E$ and any $v \in F$. From this it follows immediately that any element $u$ of $U$ can be written as $u = yu'$ with $y \in U_F = \langle y_v \mid v \in F \rangle$ and $u' \in U_E$. Clearly, $U_E$ is finite abelian and thus it is enough to show that $U_F$ is also finite in order to prove the finiteness of $U$.

We now recall that the elements $y_v$ originally came from automorphisms $g_v$ of $G_0 = \langle x_\tau \mid \tau \in T \rangle$ defined by

$$x_\tau g_v = x_\tau^{y_v}$$

(see Construction 1.1 for comparison). So, we may identify $U_F$ with $W_0 = \langle g_v \mid v \in F \rangle \subseteq \operatorname{Aut}(G_0)$. Let

$$E' = \{v \restriction n \mid n \leq k, \ v \in F\}$$

and decompose $G_0$ into $G_0 = C \times C_0$ with

$$C = \langle x_\tau \mid \tau \in E' \rangle$$

and

$$C_0 = \langle x_\tau \mid \tau \in T \setminus E' \rangle.$$

We define automorphisms $f_v, h_v$ of $G_0$ ($v \in F$) by

$$f_v \restriction C = g_v \restriction C, \ f_v \restriction C_0 = \operatorname{id}_{C_0}$$

and

$$h_v \restriction C = \operatorname{id}_C, \ h_v \restriction C_0 = g_v \restriction C_0.$$

Obviously, we thus have $g_v = f_v h_v = h_v f_v$ for all $v \in F$. Hence, if we show that

$$W = \langle f_v, h_v \mid v \in F \rangle$$

is finite, then we also deduce that $W_0 \subseteq W$ is finite. So, let us show that $W$ is finite. For any $v, w \in F$ we have $f_v h_w = h_w f_v$ by definition. Moreover, by the choice of $k$, we have that $h_v h_w = h_w h_v$ for any $v, w \in F$. Therefore, any element of $W$ can be written as $f h_{v_1} \ldots h_{v_l}$ with $v_1, \ldots, v_l \in F$ and $f \in \langle f_v \mid v \in F \rangle = W^*$. But $W^*$ can be viewed as a subgroup of $\operatorname{Aut}(C)$ and $C$ is finite. Therefore $W^*$ and thus $U$ is finite and obviously a 2-group, i.e. $H$ is a locally finite 2-group.

Finally, suppose that there is an embedding $h\colon H \hookrightarrow G$. By property $(+)$ there exists $v \in {}^\omega\lambda$ such that $h_v \subseteq h \restriction G_0 \subseteq h$. This implies that $h \restriction \langle K_v, y_v = y^v_{\varepsilon_v} \rangle$ is an embedding, contradicting the choice of $\varepsilon_v$. Therefore $H$ cannot be embedded into $G$ and so the conclusion of the theorem follows. $\square$

We finish this section with proving the following proposition using similar arguments as above.

**Proposition 1.4** *Let $\lambda$ be a cardinal satisfying $\lambda = \lambda^{\aleph_0}$. Then there exists no universal torsion group of cardinality $\lambda$, not even universal for locally finite groups.*

**Proof.** It is clear that all we need to show is the failure of the amalgamation property; the non-existence then follows as in the proof of Theorem 1.4.

To do so we again consider the locally finite 2-groups $A, B_0, B_1$ as given by Construction 1.1 (see also Proposition 1.2).

Suppose there are a group $G$ and embeddings $h_0\colon B_0 \hookrightarrow G$, $h_1\colon B_1 \hookrightarrow G$ such that $h_0 \restriction A = h_1 \restriction A$. We show that $G$ contains a torsion-free element, i.e. $G$ cannot be torsion. For any $n < \omega$ let $z_n = x_n h_0 = x_n h_1$ and $v_0 = y_0 h_0$, $v_1 = y_1 h_1$. Moreover, we put $p_{2n} = (v_1 v_0)^n$, $p'_{2n} = (v_0 v_1)^n$ and $p_{2n+1} = v_0 p_{2n} = p'_{2n} v_0 = p'_{2n+1}$. We inductively prove that, for any $k \geq 1$, the following holds:

$$p_k z_0 p'_k = z_0^{\varepsilon_0} z_1^{\varepsilon_1} \ldots z_{k-1}^{\varepsilon_{k-1}} z_k \qquad (*)$$

for some $\varepsilon_0, \varepsilon_1, \ldots, \varepsilon_{k-1} \in \{0,1\}$. First note that

$$v_0 z_i v_0 = (x_i^{y_0}) h_0 = \begin{cases} z_i z_{i+1} & ;\ \text{for } i \text{ even} \\ z_i & ;\ \text{for } i \text{ odd} \end{cases} \qquad (1)$$

and

$$v_1 z_i v_1 = (x_i^{y_1}) h_1 = \begin{cases} z_i & ;\ \text{for } i \text{ even} \\ z_i z_{i+1} & ;\ \text{for } i \text{ odd.} \end{cases} \qquad (2)$$

Now, let us consider $(*)$. For $k = 1$ we have $p_1 z_0 p'_1 = v_0 z_0 v_0 = z_0 z_1$ by (1), i.e. $(*)$ is true for $k = 1$. Now suppose that $(*)$ holds for $k \geq 1$ and consider $k+1$. Obviously, we have to distinguish between $k$ even or odd.

Let us first assume that $k$ is even, say $k = 2n$. Then, by the induction hypothesis and by (1), we have:

$$\begin{aligned} p_{2n+1} z_0 p'_{2n+1} &= v_0 p_{2n} z_0 p'_{2n} v_0 \\ &= v_0 z_0^{\varepsilon_0} \ldots z_{k-1}^{\varepsilon_{k-1}} z_k v_0 \\ &= \left( v_0 z_0^{\varepsilon_0} v_0 v_0 z_1^{\varepsilon_1} \ldots z_{k-1}^{\varepsilon_{k-1}} v_0 \right) v_0 z_k v_0 \\ &= \left( z_0^{\varepsilon'_0} z_1^{\varepsilon'_1} \ldots z_{k-1}^{\varepsilon'_{k-1}} \right) z_k z_{k+1} \end{aligned}$$

for some $\varepsilon'_i \in \{0, 1\}$ $(i = 0, \ldots k - 1)$.

Now suppose that $k$ is odd, say $k = 2n - 1$. Then, again by the induction hypothesis and by (2), we have:

$$\begin{aligned} p_{2n} z_0 p'_{2n} &= v_1 p_{2n-1} z_0 p_{2n-1} v_1 \\ &= \left( v_1 z_0^{\varepsilon_0} v_1 v_1 z_1^{\varepsilon_1} \ldots z_{k-1}^{\varepsilon_{k-1}} v_1 \right) v_1 z_k v_1 \\ &= \left( z_0^{\varepsilon'_0} z_1^{\varepsilon'_1} \ldots z_{k-1}^{\varepsilon'_{k-1}} \right) z_k z_{k+1} \end{aligned}$$

for some $\varepsilon'_i \in \{0, 1\}$ $(i = 0, \ldots k - 1)$. Therefore property $(*)$ is proven.

Finally it follows that the element $v_0 v_1 \in G$ is torsion-free since otherwise

$$p'_{2n} = (v_0 v_1)^n = 1 = (v_1 v_0)^n = p_{2n}$$

for some $n \geq 1$ and so

$$x_0 h_0 = z_0 = p_{2n} z_0 p'_{2n} = z_0^{\varepsilon_0} z_1^{\varepsilon_1} \ldots z_{2n-1}^{\varepsilon_{2n-1}} z_{2n} = \left( x_0^{\varepsilon_0} x_1^{\varepsilon_1} \ldots x_{2n-1}^{\varepsilon_{2n-1}} \right) h_0 \; x_{2n} h_0$$

contradicting the injectivity of $h_0$ together with the unique representations of elements of $A$ as product of the $x_i$'s. □

We would like to mention that the methods used in this section follow the pattern of [GrS]. In particular, the result that there is no universal locally finite group of cardinality $\lambda = \lambda^{\aleph_0}$ follows from Proposition 1.4 (see [GrS, Theorem 12]).

## §2 Existence of epi-universal locally nilpotent groups

Here we show that, for certain cardinals $\lambda$, there exist epi-universal groups in $\mathcal{LN}_\lambda$. In fact, we shall construct an epi-universal group which is also torsion-free. As mentioned before we thus obtain an affirmative answer to a problem by Plotkin (see [KN, 3.47]).

First we define locally nilpotent groups which are freely generated except the necessary relations.

**Definition 2.1** *Let $M$ be any set and $c: [M]^{<\aleph_0} \to \omega$ an increasing function where $[M]^{<\aleph_0}$ denotes the family of all finite subsets of $M$. For all $U \subseteq M$ we put $F_U = \langle x_m \mid m \in U \rangle$ to be the free group generated by the $x_m$'s; we write $F = F_M$. We define the* pseudo-free locally nilpotent group $\mathrm{Fr}(M, c)$ *by*

$$\mathrm{Fr}(M, c) = F/K$$

*where $K = K(M, c)$ is the normal subgroup generated by*

$$\{F_U^{Uc} \mid U \in [M]^{<\aleph_0}\}.$$

Recall, for any group $G$, $G^n$ is inductively defined by $G^0 = G$, $G^{n+1} = [G^n, G]$ and $G^n$ is generated by $\{[g_0, \ldots, g_n] \mid g_i \in G\}$ where the commutator brackets accumulate on the left.

It is immediate from the definition that $\mathrm{Fr}(M, c)$ is, indeed, locally nilpotent.

We order the pairs $(M, c)$ by $(M_1, c_1) \leq (M_2, c_2)$ if $M_1 \subseteq M_2$ and $c_2$ extends $c_1$. For such ordered pairs we first show:

**Lemma 2.2** *Let $(M_1, c_1) \leq (M_2, c_2)$ be as above. Then*
$$\mathrm{Fr}(M_2, c_2) = D \rtimes G_1 \text{ where } G_1 \cong \mathrm{Fr}(M_1, c_1).$$

**Proof.** Let $F_1 = F_{M_1}, F_2 = F_{M_2}, K_1 = K(M_1, c_1), K_2 = K(M_2, c_2), H_1 = \mathrm{Fr}(M_1, c_1) = F_1/K_1$ and $H_2 = \mathrm{Fr}(M_2, c_2) = F_2/K_2$.

We define $h\colon F_2 \to F_1$ to be the canonical projection, i.e.

$$x_m h = \begin{cases} x_m, & \text{for } m \in M_1 \\ 1, & \text{else.} \end{cases}$$

Obviously, $K_1 \subseteq K_2 h$ by Definition 2.1. Conversely, a typical generator of $K_2$ looks like $k^f = f^{-1}kf$ for some $f \in F_2$ and $k = [k_0, \ldots, k_n]$ with $k_i \in F_U$, $i \leq n = Uc_2$ for some finite $U \subseteq M_2$. Then $kh = [k_0, \ldots, k_n]h = [k_0h, \ldots, k_nh]$ where $k_i h \in F_{U \cap M_1}$ and $(U \cap M_1)c_1 = (U \cap M_1)c_2 \leq Uc_2 = n$. Hence $kh \in F_{U \cap M_1}^n \subseteq F_{U \cap M_1}^{(U \cap M_1)c_1} \subseteq K_1$ and so $(k^f)h = (kh)^{fh} \in K_1$ since $fh \in F_1$.

Therefore we have shown that $K_1 = K_2 h$ and thus $h$ induces an epimorphism
$$\overline{h}\colon F_2/K_2 = H_2 \to H_1 = F_1/K_1.$$
Moreover, $h \restriction F_1 = \mathrm{id}_{F_1}$ and $K_2 h = K_1$ imply $K_1 = F_1 \cap K_2$ and thus $H_1 \cong (F_1 K_2)/K_2 =: G_1$. Also, $G_1 \cap \operatorname{Ker} \overline{h} = 1$ since, for $f \in F_1$ with $fK_2 \in \operatorname{Ker} \overline{h} = (\operatorname{Ker} h\ K_2)/K_2$, we obtain $fK_2 = fhK_2 = K_2$. Hence $H_2 = F_2/K_2 = \operatorname{Ker} \overline{h} \rtimes G_1$ since $F_2 = F_1 * F_{M_2 \setminus M_1}$ and $F_{M_2 \setminus M_1} \subseteq \operatorname{Ker} h$. Thus the assertion follows with $D = \operatorname{Ker} \overline{h}$. □

Next we show that the pseudo-free locally nilpotent groups are torsion-free. For this we need the commutator collecting process due to P. Hall (see [H, Ch. 10]). For the convenience of the reader we here recall the definition of basic commutators and the main result concerning them. However, we shall not go into the details of the "collecting process".

**Definition 2.3** *Let $F$ be the free group generated by $x_0, x_1, \ldots, x_r$ $(r < \omega)$. We define the* basic commutators *(for $F$) of weight $n \geq 1$ as follows:*

**(1)** *The $x_i$'s $(i = 0, 1, \ldots, r)$ are the basic commutators of weight 1;*

**(2)** *Having defined the basic commutators of weight less than $n$, the basic commutators of weight $n$ are the elements of the form $[y_1, y_2]$ where*

   (a) *$y_1, y_2$ are basic commutators of weights $n_1, n_2$, respectively, with $n_1 + n_2 = n$, and*

   (b) *$y_1 > y_2$, and if $y_1 = [z_1, z_2]$ with $z_1, z_2$ basic commutators, then $y_2 \geq z_2$*

*where the commutators are ordered according to their weight, but arbitrarily well ordered among those of the same weight.*

The use of the word "basic" in the above definition is due to the following result. For a proof we refer to [H, Theorem 11.2.4].

**Lemma 2.4** *Let $F$ be a finitely generated free group and let $b_0 < b_1 < \ldots < b_t$ be the basic commutators of weight $\leq n + 1$.*

*Then, any element $f$ of $F$ has a unique representation*

$$f \equiv b_0^{e_0} \ldots b_t^{e_t} \mod F^{n+1}$$

*($e_0, \ldots e_t \in \mathbb{Z}$). Moreover, the basic commutators of weight $n+1$ form a basis for the free abelian group $F^n/F^{n+1}$.* □

We are now ready to prove the torsion-freeness of the pseudo-free locally nilpotent groups.

**Proposition 2.5** *Let $(M, c)$ and $\mathrm{Fr}(M, c)$ be as in Definition 2.1. Then $\mathrm{Fr}(M, c)$ is a torsion-free group.*

**Proof.** By Lemma 2.2 it is enough to consider $\mathrm{Fr}(M, c) = F_M/K(M, c)$ for a finite set $M$. So, let $M$ be finite, $c \colon \mathcal{P}(M) \to \omega$ increasing, $F = F_M$ and $K = K(M, c)$ where $\mathcal{P}(M)$ denotes the power set of $M$. Then $K$ is the normal subgroup generated by the set $\{F_U^{Uc} \mid U \subseteq M\}$, in particular we have that $F^{Mc} \subseteq K$.

First we show that, for any $n < \omega$, the group

$$\overline{K}_n = ((K \cap F^n) F^{n+1})/F^{n+1} \text{ is a direct summand} \tag{+}$$
$$\text{of the free abelian group } F^n/F^{n+1}.$$

To do so we define basic commutators for $K$ as follows (compare Definition 2.3):

(1) $x_m$ is a *basic commutator for $K$ of weight* 1 if $\{m\}c = 0$.

(2) Having defined the basic commutators for $K$ of weight less than $n$ ($\geq 1$), the *basic commutators for $K$ of weight $n$* are the basic commutators for $F_U$ of weight $n$ for all $U \subseteq M$ with $Uc < n$ and the elements of the form $[y_1, y_2]$ such that

(a) $y_1, y_2$ are basic commutators for $K$ of weights $n_1, n_2$ with $n_1 + n_2 = n$, and

(b) $y_1 > y_2$, and if $y_1 = [z_1, z_2]$ with $z_1, z_2$ basic commutators for $K$, then $y_2 \geq z_2$.

It follows immediately from this definition that the basic commutators for $K$ form a subset of the basic commutators (for F); we may assume that the ordering of the basic commutators for the $F_U$'s, for $K$ and for $F = F_M$ is the same.

We show that the basic commutators for $K$ of weight $n+1$ form a basis for the free abelian group $\overline{K}_n$; the assertion (+) then follows.

Clearly, the basic commutators for $K$ of weight $n + 1$ are linearly independent elements of $\overline{K}_n$. It remains to show that they generate $\overline{K}_n$.

For any $k \in F^n$ and any $f \in F$ we have $k^f \equiv k[k, f] \equiv k \mod F^{n+1}$ since $[k, f] \in F^{n+1}$. Hence

$$\overline{K}_n = \left(\left(\left\langle F_U^{Uc} \mid U \subseteq M,\ Uc \leq n\right\rangle \cap F^n\right) F^{n+1}\right) / F^{n+1}.$$

It follows from Lemma 2.4 that any element of $F_U^{Uc}$ can be uniquely expressed as a product of powers of basic commutators for $F_U$ of weights $Uc+1, \ldots, n+1$ modulo $F_U^{n+1} \subseteq F^{n+1}$. By definition these basic commutators are also basic commutators for $K$. Therefore, any element $f$ of $\left\langle F_U^{Uc} \mid U \subseteq M,\ Uc \leq n\right\rangle$ can be expressed as a product of powers of basic commutators for $K$ of weights $\leq n+1$ modulo $F^{n+1}$. To obtain the correct ordering, and thus a unique representation, we replace a product $cd$ with $c > d$ by $dc[c, d]$ where $[c, d]$ is also a basic commutator for $K$. On the other hand, if $f \in F^n$ then $f$ can be uniquely expressed as a product of powers of basic commutators for $F$ of weight $n+1$ modulo $F^{n+1}$. Hence, for an element $f$ of $\left\langle F_U^{Uc} \mid U \subseteq M,\ Uc \leq n\right\rangle \cap F^n$ we conclude, comparing both possible representations, that $f$ is a product of powers of basic commutators for $K$ of weight $n+1$, i.e. $\overline{K}_n$ is generated by these commutators and so (+) follows.

Finally we show, by induction on $n$, that $F/(KF^n)$ is torsion-free for any $n < \omega$. It then follows that $\mathrm{Fr}(M, c) = F/K$ is torsion-free since $KF^n = K$ for $n \geq Mc$. For $n = 0$ we have $F/(KF^0) = F/F = 1$ which is trivially torsion-free.

Suppose $F/(KF^n)$ is torsion-free and consider the canonical epimorphism $\pi : F/(KF^{n+1}) \to F/(KF^n)$. Then the torsion-freeness of $F/(KF^{n+1})$ follows from the torsion-freeness of $\mathrm{Ker}\,\pi$. Now,

$$\mathrm{Ker}\,\pi = (KF^n)/(KF^{n+1}) \cong F^n/((K \cap F^n)F^{n+1})$$

where the latter is isomorphic to
$$(F^n/F^{n+1})/(((K \cap F^n)F^{n+1})/F^{n+1}) = (F^n/F^{n+1})/\overline{K}_n$$
which is torsion-free, in fact, free abelian by (+). Therefore $F/(KF^n)$ is torsion-free for any $n$ and hence so is $\text{Fr}(M, c) = F/K$ as required. □

Next we show that the pseudo-free $\mathcal{LN}$-groups have some property in common with $\mathcal{V}$-free groups in a variety $\mathcal{V}$.

**Proposition 2.6** *Any locally nilpotent group is an epimorphic image of a pseudo-free locally nilpotent group.*

**Proof.** Let $G \in \mathcal{LN}$ and let $\{y_i \mid i \in I\}$ be a set of generators for $G$. Moreover, let $Uc$ be the nilpotency class of $G_U = \langle y_i \mid i \in U \rangle$ for any finite $U \subseteq I$. Then $c_G = c \colon [I]^{<\aleph_0} \to \omega$ is an increasing function.

We show that $G$ is an epimorphic image of the pseudo-free locally nilpotent group $\text{Fr}(I, c) = F/K$ (see Definition 2.1). First we define $\pi \colon F = F_I \to G$ by $x_i\pi = y_i$ ($i \in I$). Clearly, $\pi$ is surjective. Moreover, a typical generator of $K$ is of the form $k^f$ for $f \in F$ and $k = [k_0, \ldots, k_n]$ with $k_i \in F_U = \langle x_i \mid i \in U \rangle$ for some finite $U \subseteq I$ with $Uc = n$. Hence $k\pi = [k_0\pi, \ldots, k_n\pi]$ where $k_i\pi \in F_U\pi = G_U$ and thus $k\pi = 1 = (k\pi)^{(f\pi)} = (k^f)\pi$ since $Uc$ is the nilpotency class of $G_U$. Therefore $K\pi = 1$ and so $\pi$ induces an epimorphism $\overline{\pi} \colon \text{Fr}(I, c_G) = F/K \to G$ which completes the proof. □

As an immediate consequence from the Propositions 2.5 and 2.6 we obtain an affirmative answer to Plotkin's question:

**Corollary 2.7** *Any locally nilpotent group is an epimorphic image of a torsion-free locally nilpotent group.* □

We are finally ready to prove the main result of this section:

**Theorem 2.8** *Let $\lambda$ be a cardinal such that $\lambda = \lambda^{<\lambda}$. Then there exists an epi-universal group in the class $\mathcal{LN}_\lambda$.*

Note, assuming the generalized continuum hypothesis GCH, any regular cardinal $\lambda$ satisfies $\lambda = \lambda^{<\lambda}$.

**Proof of 2.8:** Let $\lambda$ be as above. First we consider the class

$$\mathcal{K} = \{(M, c) \mid \emptyset \neq M \text{ a set}, c \colon [M]^{<\aleph_0} \to \omega \text{ an increasing function}\}.$$

It is easy to check that the class $\mathcal{K}$ satisfies the assumptions of Jónsson's Theorem (see [BS, p.213]). Hence there exists an universal element $(M_u, c_u) \in \mathcal{K}_\lambda$ where $\mathcal{K}_\lambda$ is the subclass of $\mathcal{K}$ consisting of those elements $(M, c)$ with $|M| \leq \lambda$.

We define $G_u$ to be the pseudo-free locally nilpotent group $G_u = \mathrm{Fr}(M_u, c_u)$ (see Definition 2.1). Then $G_u$ is clearly an element of $\mathcal{LN}_\lambda$. We show that $G_u$ is epi-universal in $\mathcal{LN}_\lambda$. To do so let $G \in \mathcal{LN}_\lambda$ be arbitrary. By Proposition 2.6 there are a pseudo-free locally nilpotent group $\widetilde{G} = \mathrm{Fr}(I, c_G)$ and an epimorphism $\pi \colon \widetilde{G} \twoheadrightarrow G$ where $|\widetilde{G}| = |I| \cdot \aleph_0 = |G| \leq \lambda$, i.e. $\widetilde{G} \in \mathcal{LN}_\lambda$ and $(I, c_G) \in \mathcal{K}_\lambda$.

Now, since $(M_u, c_u)$ is universal in $\mathcal{K}_\lambda$, there exists $M \subseteq M_u$ such that $(M, c = c_u \upharpoonright [M]^{<\aleph_0})$ is isomorphic to $(I, c_G)$, i.e. there is a bijection $\varphi \colon M \to I$ such that $(Nc)\varphi = (N\varphi)c_G$ for any finite subset $N$ of $M$. Therefore $\varphi$ induces an isomorphism $\overline{\varphi} \colon F_M \to F_I$ such that $K(M, c)\overline{\varphi} = K(I, c_G)$ by the definition of the $K$'s (see Definition 2.1). Hence $\overline{\varphi}$ induces $\widetilde{\varphi} \colon \mathrm{Fr}(M, c) = F_M/K(M, c) \to F_I/K(I, c_G) = \mathrm{Fr}(I, c_G) = \widetilde{G}$ which is obviously also an isomorphism.

Finally, there is an epimorphism $\pi' \colon G_u \to \mathrm{Fr}(M, c)$ by Lemma 2.2 and thus the composition $\pi'\widetilde{\varphi}\pi$ is an epimorphism from $G_u$ onto $G$, i.e. $G_u$ is epi-universal as required. $\square$

## §3 An observation

In this final section we slightly generalize the methods from §2 to other "local classes". We observe that results similar to Lemma 2.2, Proposition 2.6 and Theorem 2.8 hold. However, since the proof of the torsion-freeness of the pseudo-free locally nilpotent groups heavily uses typical properties of nilpotent groups (see Proposition 2.5), there will be no corresponding statement in this more general context.

Throughout let $\{\mathcal{V}_n \mid n < \omega\}$ be an ascending chain of group varieties. Moreover, let $E_n$ be a set of identities generating the variety $\mathcal{V}_n$ ($n < \omega$) where we assume $E_n \supseteq E_{n+1}$. We put $\mathcal{V} = \bigcup_{n<\omega} \mathcal{V}_n$ and define $\mathcal{LV}$ to be the class of all $\mathcal{V}$-*locally* groups (or $\mathcal{LV}$-*groups*), that is, all groups $G$ such that any finitely generated subgroup of $G$ belongs to $\mathcal{V}$.

As before in §2 we first define pseudo-free $\mathcal{LV}$-groups:

**Definition 3.1** *Let $\{E_n \mid n < \omega\}$ be the descending chain of the sets $E_n$ of identities, say $E_n = \{\tau = 1 \mid \tau = \tau(X_0, \ldots, X_{n_\tau}) \in T_n\}$. Moreover, let $M$ be any set and let $c : [M]^{<\aleph_0} \to \omega$ be an increasing function. For all $U \subseteq M$ let $F_U = \langle x_m \mid m \in U \rangle$ and put $F = F_M$.*

*We define the pseudo-free $\mathcal{LV}$-group $\mathrm{Fr}_{\mathcal{LV}}(M, c)$ by*

$$\mathrm{Fr}_{\mathcal{LV}}(M, c) = F/K$$

*where $K = K_{\mathcal{LV}}(M, c)$ is the normal subgroup generated by*

$$\bigcup_{U \in [M]^{<\aleph_0}} \{\tau(f_0, \ldots, f_{n_\tau}) \mid \tau \in T_{Uc}, \ f_0, \ldots, f_{n_\tau} \in F_U\}.$$

It follows immediately from the definition that $\mathrm{Fr}_{\mathcal{LV}}(M, c)$ is an $\mathcal{LV}$-group since, for any finite set $U \subseteq M$, we have that $\tau(\bar{f}_0, \ldots, \bar{f}_{n_\tau}) = 1$ for any $\tau \in T_{Uc}$ and $\bar{f}_0, \ldots, \bar{f}_{n_\tau} \in (F_U K)/K$, i.e. $(F_U K)/K \in \mathcal{V}_{Uc}$.

We order the pairs $(M, c)$ as before.

**Lemma 3.2** *If $(M_1, c_1) \leq (M_2, c_2)$ then*

$$\mathrm{Fr}_{\mathcal{LV}}(M_2, c_2) = D \rtimes G_1 \text{ where } G_1 \cong \mathrm{Fr}_{\mathcal{LV}}(M_1, c_1).$$

**Proof.** The proof is the same as the one of Lemma 2.2 replacing the commutators by the identities $\tau = 1 \in E_n$ using the fact that $E_n \supseteq E_{n+1}$. $\square$

The next result is the corresponding one to Proposition 2.6. Although the arguments to prove it are basically the same as before, we want to include the proof to illustrate the minor differences.

**Proposition 3.3** *Any $\mathcal{LV}$-group is an epimorphic image of some pseudo-free $\mathcal{LV}$-group.*

**Proof.** Let $G$ be any $\mathcal{LV}$-group and let $\{y_i \mid i \in I\}$ be a set of generators for $G$. Moreover, for any finite $U \subseteq I$, let $Uc$ be the minimal $n < \omega$ such that $\langle y_i \mid i \in U \rangle \in \mathcal{V}_n$, which exists since $G \in \mathcal{LV}$, $\mathcal{V} = \bigcup_{n<\omega} \mathcal{V}_n$. Then $c_G = c\colon [I]^{<\aleph_0} \to \omega$ is increasing since any $\mathcal{V}_n$ is closed under subgroups.

We show that $G$ is an epimorphic image of the pseudo-free $\mathcal{LV}$-group $\mathrm{Fr}_{\mathcal{LV}}(I,c)$. As in the proof of Proposition 2.6 we first define an epimorphism $\pi\colon F = F_I \to G$ by $x_i\pi = y_i$. Now, a typical generator of $K = K_{\mathcal{LV}}(I,c)$ is of the form $k^f$ for $f \in F$ and $k = \tau(k_0, \ldots, k_{n_\tau})$ with $\tau \in T_{Uc}$, $k_0, \ldots, k_{n_\tau} \in F_U$ for some finite $U \subseteq I$. Then $k\pi = \tau(k_0\pi, \ldots, k_{n_\tau}\pi) = 1 = (k^f)\pi$ since $k_i\pi \in \langle y_i \mid i \in U \rangle \in \mathcal{V}_n$. Therefore $K\pi = 1$ and so $\pi$ induces an epimorphism $\pi'\colon \mathrm{Fr}_{\mathcal{LV}}(I,c_G) = F/K \to G$ which completes the proof. □

Finally, the existence of epi-universal $\mathcal{LV}$-groups follows by exactly the same arguments as in the proof of Theorem 2.8 using Jónsson's Theorem, Lemma 3.2 and Proposition 3.3.

**Theorem 3.4** *Let $\lambda$ be a cardinal such that $\lambda = \lambda^{<\lambda}$. Then there exists an epi-universal group in the class $\mathcal{LV}_\lambda$ of all $\mathcal{LV}$-groups of cardinality $\leq \lambda$.* □

Notice that, as mentioned in the introduction, we thus obtain the existence of, for example, epi-universal locally solvable and epi-universal torsion groups. Further examples are left to the reader.

Rüdiger Göbel, Simone L. Wallutis                Saharon Shelah
FB 6, Mathematik und Informatik                  Department of Mathematics
Universität Essen                                Hebrew University
45117 Essen, Germany                             Jerusalem, Israel
Email: r.goebel@uni-essen.de                     Email: shelah@math.huji.ac.il
       simone.wallutis@uni-essen.de